\documentclass[12pt,oneside]{article}
\usepackage{amsmath,amssymb,amsfonts,amsthm}
\usepackage{color}
\newcommand{\set}[1]{\left\{#1\right\}}

\def\Section#1{\vspace{30truept}\addtocounter{section}{1}\setcounter{thm}{0}
\setcounter{equation}{0}{\noindent\Large\bf
    \arabic{section}.~~#1}\par \vspace{12pt}}

\newtheorem{thm}{Theorem}[section]
\newtheorem{cor}[thm]{Corollary}
\newtheorem{lem}[thm]{Lemma}
\newtheorem{prop}[thm]{Proposition}
\newtheorem{defn}[thm]{Definition}

\newtheorem{rem}[thm]{Remark}

\numberwithin{equation}{section}
\newcommand\overcirc[1]{\raisebox{10pt}{\tiny{$\circ$}}{\kern-7.5pt}\mbox{$#1$}}
\newcommand\undersym[2]{\raisebox{-6pt}{$#2$}{\kern-5pt}\mbox{$#1$}}
\newcommand\overdiamond[1]{\raisebox{10pt}{\small$\star$}{\kern-7.5pt}\mbox{$#1$}}
\newcommand\overast[1]{\raisebox{10pt}{\small$\ast$}{\kern-7.5pt}\mbox{$#1$}}
\newcommand\overlind[1]{\raisebox{10pt}{\small$\overline{{\hspace{2pt}}\star}$}{\kern-7.5pt}\mbox{$#1$}}
\newcommand\overlinc[1]{\raisebox{10pt}{\small$\overline{{\hspace{2pt}}\circ}$}{\kern-7.5pt}\mbox{$#1$}}
\newcommand\overlina[1]{\raisebox{10pt}{\small$\overline{{\hspace{1pt}}\ast}$}{\kern-7.5pt}\mbox{$#1$}}

\begin{document}

\title{{\bf{ Generalized $\beta$-conformal change and special Finsler spaces}}}
%\footnote{ArXiv Number: 1004.5478 [math.DG]}
\author{{\bf Nabil L. Youssef$^\dag$, S. H. Abed$^\ddag$
 and S. G. Elgendi$^\sharp$}}
\date{}
%\thanks{\it Department of Mathematics, etc}
%\pagestyle{fancy}

             % End of preamble and beginning of text.
\maketitle                     % Produces the title.
\vspace{-1.15cm}
\begin{center}
{$^\dag$Department of Mathematics, Faculty of Science,\\
Cairo University, Giza, Egypt}
\end{center}
\vspace{-0.8cm}
\begin{center}
 nlyoussef2003@yahoo.fr, nyoussef@frcu.eun.eg
\end{center}

\begin{center}
{$^\ddag$Department of Mathematics, Faculty of Science,\\
Cairo University, Giza, Egypt}
\end{center}
\vspace{-.8cm}
\begin{center}
sabed52@yahoo.fr, sabed@frcu.eun.eg
\end{center}
\vspace{-.3cm}
\begin{center}
{$^\sharp$Department of Mathematics, Faculty of Science,\\
Benha University, Benha,
 Egypt}
\end{center}
\vspace{-.8cm}
\begin{center}
salahelgendi@yahoo.com
\end{center}
\smallskip

\vspace{1cm} \maketitle
\smallskip

\noindent{\bf Abstract.} In this paper, we  investigate  the change
of Finslr metrics  $$L(x,y) \longrightarrow\overline{L}(x,y) =
f(e^{\sigma(x)}L(x,y),\beta(x,y)),$$ which we   refer to as a
generalized $\beta$-conformal change.  Under this change, we study
some special Finsler spaces, namely, quasi C-reducible, semi
C-reducible, C-reducible,  $C_2$-like, $S_3$-like and $S_4$-like
Finsler spaces. We obtain some characterizations of the energy
$\beta$-change, the Randers change and the Kropina change. We also obtain the
transformation of the  T-tensor under this change and study some
interesting special cases. We then impose a certain condition on the
generalized $\beta$-conformal change, which we call the b-condition,
and investigate the geometric consequences of such a condition.
Finally, we give the conditions under which a generalized
$\beta$-conformal change is projective and generalize  some known
results in the literature.

\newpage

%%%%%%%%%%%%%%%%%%%%%%%%%%%%%%%%%%%%%%%%%%%%%%%%%%%%%%%%%%% Introduction %%%%%%%%%%%%%%%%%%%%%%%%%%%%%%%%%%%%%%%%%%%%%%%%%%%%%%%%%%%

\vspace{30truept}\centerline{\Large\bf{Introduction}}\vspace{12pt}
\par
 Throughout,  $M$ is  an n-dimensional
 $C^\infty$ differentiable
manifold and $F^n=(M,L)$  a Finsler space equipped with the
fundamental function $L(x,y)$ on $TM$.
\par Finsler
geometry was first introduced by Finsler himself, to be studied by
many eminent mathematicians for its theoretical importance and
applications in the variational calculus, mechanics and theoretical
physics. Moreover, the dependence of the fundamental function $L(x,
y)$ on both the positional argument $x$ and directional argument $y$
offers the possibility to use it to describe the {anisotropic}
properties of the physical space. For a differential one-form
$\beta(x,dx)=b_{i}(x)dx^{i}$ on $ M$, Randers \cite{r2.11}, in 1941,
introduced a special Finsler space with the Finsler change
$\overline{L}=L+\beta$, where $L$ is Riemannian, to consider a
unified field theory ($L=\sqrt{a_{ij}y^iy^j}$, $a_{ij}$
  being   the gravitational tensor
field and $b_{i}(x)$ the electromagnetic potential). Masumoto
\cite{r2.8}, in 1974, studied Randers  and generalized Randers
changes in which $L$ is Finslerian. Kropina \cite{r1.4} introduced
the change $\overline{L}=L^2/\beta$, where $L$ is Riemannian. This
change has been studied by other authors such as Shibata
\cite{r2.12} and Matsumoto \cite{r2.9}. Randers  and Kropina changes
are closely related to physical theories and so Finsler spaces with
these metrics have been investigated by many authors, from various
approaches in both the physical and mathematical aspects
(\cite{r2.5}, \cite{r2.8},\cite{r2.14}, \cite{r2.16},  \cite{Tay}, \cite{r2.17},
\cite{r2.18}). Randers change was also applied to the theory of the
electron microscope by R. S. Ingarden \cite{Ingarden}. Moreover,
there is some relation between the Kropina metric and the Lagrangian
function of analytic dynamics \cite{r2.12}.  In 1984, Shibata
\cite{r2.13} considered  the general case of any $\beta$-change,
that is, $\overline{L}=f(L,\beta)$, thus generalizing   many changes
in Finsler geometry  (\cite{r2.8}, \cite{r2.12}). In this context,
he studied some special Finsler spaces, such as C-reducible and
$S_4$-like, under Randers change.

 On the other hand, in 1976,
Hashiguchi \cite{r2.6} studied  the conformal change of a Finsler
metric, namely, $ \overline{L} = e^{\sigma(x)} L$. In particular, he
also dealt with the special conformal transformation named
C-conformal. This change has been studied by Izumi \cite{r2.7} among
others. In 2008, Abed (\cite{r2.1}, \cite{r2.2}) introduced the
change $\bar L = e^{\sigma(x)} L + \beta$,   which he called a
$\beta$-conformal change, thus generalizing the conformal, Randers
and generalized Randers changes. Moreover, he studied some special
Finsler space under this change such as C-reducible, $S_3$-like and
$S_4$-like.

In \cite{gomaa}, the present authors introduced and investigated the more general
change of Finsler metrics: $$ L(x,y)
\rightarrow\overline{L}=f(e^{\sigma(x)}L(x,y),\beta(x,y))=f(\widetilde{L},\beta),
$$ where  $\widetilde{L}= e^{\sigma(x)}L$ and $f$ is a positively
homogeneous function of $\widetilde{L}$ and $\beta$ of degree one.
They obtained the difference between Cartan connection associated
with $(M,L)$ and Cartan connection associated with
$(M,\overline{L})$, also, they established some interesting results
and computed the torsion and curvature tensors of the transformed
space $(M,\overline{L})$ for the four fundamental connections in
Finsler geometry. This change is referred to as a generalized
$\beta$-conformal change. It is clear that this change is a
generalization of all the above mentioned changes and deals
simultaneously with $\beta$-change and conformal change. It combines both
cases of Shibata ($\overline{L}=f(L,\beta)$) and that of
Hashiguchi ($\overline{L}=e^\sigma L$).

In this paper, we  continue our investigation of the generalized
$\beta$-conformal change.
Under this change, we study some special Finsler spaces, compute the
transformed T-tensor, introduce what we call b-condition and study when this
change becomes projective.

The present paper is organized as follows. In section 1, we
introduce the necessary material and background required for the
present work. In section 2, we deal with some special Finsler spaces
under a generalized $\beta$-conformal change, namely, quasi
C-reducible, Semi C-reducible, C-reducible, $C_2$-like, $S_3$-like
and $S_4$-like. In section 3, we compute the T-tensor of the
transformed space under a generalized $\beta$-conformal change and
 study some interesting special cases. In section 4, we
impose a certain condition on the generalized $\beta$-conformal change, which we call
the b-condition, and investigate the geometric consequences of such a condition.
Finally, in section 5, we give the conditions
under which a generalized $\beta$-conformal change is
projective and  generalize  some known results in the
literature.

%%%%%%%%%%%%%%%%%%%%%%%%%%%%%%%%%%%%%%%%%%%%%%%%%%% 1. Notations and Prelemenaries %%%%%%%%%%%%%%%%%%%%%%%%%%%%%%%%%%%%%%%%%%%%%%%%%%%%%%%%%%%

\Section{\textbf{  Notations and preliminaries}}
 Throughout the
present paper  we use  the terminology and notations of
\cite{gomaa}. Let  $(M,L)$ be  an n-dimensional $C^\infty$ Finsler
manifold; L being the fundamental Finsler function. Let $(x^i)$ be
the coordinates of any  point of the base manifold M  and $(y^i)$ a
supporting element at the same point.  We use the following
notations:\\
   $\partial_i$: partial differentiation with respect to $x^i$,\\
   $\dot{\partial}_i$:  partial differentiation
    with respect to  $y^i$ (basis vector fields of the vertical bundle),\\
   $g_{ij}:=\frac{1}{2}\dot{\partial}_i\dot{\partial}_j L^2=\dot{\partial}_i
   \dot{\partial}_jE$:  the
   Finsler metric tensor;  $E:=\frac{1}{2}L^2$:   the energy function,\\
   $l_i:=\dot{\partial}_iL=g_{ij} l^j=g_{ij}\frac{y^j}{L}$: the
    normalized supporting element; $l^i:=\frac{y^i}{L}$,\\
 $l_{ij}:=\dot{\partial}_il_j$,\\
   $h_{ij}:=Ll_{ij}=g_{ij}-l_il_j$:  the angular metric tensor,\\
   \vspace{7pt}$C_{ijk}:=\frac{1}{2}\dot{\partial}_kg_{ij}=\frac{1}{4}\dot{\partial}_i
    \dot{\partial}_j\dot{\partial}_k L^2$:  the Cartan  tensor, \\
   $G^i$: the components of the   canonical spray associated
    with $(M,L)$,\\
    $ N^i_j:=\dot{\partial}_jG^i$: the Barthel (or Cartan nonlinear) connection
    associated with $(M,L)$,\\
    $\delta_i:=\partial_i-N^r_i\dot{\partial}_r$: the basis
     vector fields of the horizontal bundle,\\
     \vspace{7pt}$G^i_{jh}:=\dot{\partial}_hN^i_j=\dot{\partial}_h\dot{\partial}_jG^i$:
     the coefficients of Berwald connection, \\
      $C^i_{jk}:=g^{ri}C_{rjk}=\frac{1}{2}g^{ir}\dot{\partial}_kg_{rj}$:
     the h(hv)-torsion tensor,\\
     \vspace{7pt}$\gamma^i_{jk}:=\frac{1}{2}g^{ir}(\partial_jg_{kr}+\partial_kg_{jr}-\partial_rg_{jk})$:
        the Christoffel symbols with respect to $\partial_i$,\\
        \vspace{.3cm}$\Gamma^i_{jk}:=\frac{1}{2}g^{ir}(\delta_jg_{kr}+\delta_kg_{jr}-\delta_rg_{jk})$:
          the Christoffel symbols with respect to $\delta_i$,\\
       $(\Gamma^i_{jk},N^i_j,C^i_{jk})$: The Cartan connection
        $C\Gamma$.\\

    For a Cartan  connection $(\Gamma^i_{jk},
         N^i_j,C^i_{jk})$, we define\\
         $X^i_{j\mid k}:=\delta_kX^i_j+X^m_j\Gamma^i_{mk}-X^i_m\Gamma^m_{jk}$: the
    horizontal covariant derivative of $X^i_j$,\\
   $ X^i_j|_k:=\dot{\partial}_kX^i_j+X^m_jC^i_{mk}-X^i_mC^m_{jk}$: the
    vertical  covariant derivative of $X^i_j$.\\

 Let $F^n=(M,L)$  be an n-dimensional Finsler space with a
fundamental function $L= L(x,y)$. Consider the following change of
Finsler structures,  which will be called a generalized
$\beta$-conformal change,
\begin{equation}\label{change} L(x,y)
\rightarrow\\\overline{L}(x,y)=f(e^{\sigma(x)}L(x,y),\beta(x,y))=f(\widetilde{L},\beta),
\end{equation}
where  $f$ is a positively homogeneous function of degree one in
$\widetilde{L}=e^\sigma L$ and  $\beta$  and \,\,
$\beta=b_i(x)dx^i$. Assume that $\overline{F}^n=(M,\overline{L})$
has the structure of a Finsler space. Entities related to
$\overline{F}^n$ will be denoted by barred symbols.

We define
   $$f_1:=\frac{\partial f}{\partial \widetilde{L}}\,,\qquad
   f_2:=\frac{\partial f}{\partial \beta}\,,\qquad
f_{12}:=\frac{\partial^2 f}{\partial \widetilde{L}\partial\beta},
\cdots \,\text{etc.},$$ where $\widetilde{L}=e^\sigma L$. We use
the following notations:
 \begin{eqnarray*}
   q&:=&ff_2,\hspace{3.8cm} p:=ff_1/L,\\
    q_0&:=&ff_{22},\hspace{3.5cm} p_0:=f^2_2+q_0,\\
    q_{-1}&:=&ff_{12}/L,\hspace{2.6cm} p_{-1}:= q_{-1}+(pf_2/f),\\
    q_{-2}&:=&f(e^{\sigma}f_{11}-(f_{1}/L))/L^2, \quad p_{-2}:=q_{-2}+(e^\sigma p^2/f^2).
 \end{eqnarray*}
  Note that the subscript under each of    the above geometric  objects indicates
the degree of homogeneity of that   object. We also use the
notations:
$$ m_i:= b_i-(\beta/L^2)y_i\neq 0,
\quad p_{02}:=\frac{\partial p_0}{\partial \beta}.$$
\begin{prop}\label{h-g.1} Under a generalized $\beta$-conformal change,
we have:\vspace{-0.2cm}
\begin{description}
    \item[(a)]
 ${\quad}\overline{l}_i=e^\sigma f_1 l_i+f_2 b_i$,

    \item[(b)]
${\quad}\overline{h}_{ij}=e^\sigma p \,h_{ij}+q_0 m_i m_j$,

    \item[(c)]
        ${\quad}\overline{g}_{ij}=e^\sigma p\,g_{ij}+p_0\,b_ib_j
        +e^\sigma p_{-1}(b_iy_j+b_jy_i)+
        e^\sigma p_{-2}\,y_i\,y_j.$

\item[(d)]
 The inverse metric $\overline{g}^{ij}$ of the metric
 $\overline{g}_{ij}$ is given by
    \begin{eqnarray*}
\overline{g}^{ij}
&=&(e^{-\sigma}/p){g}^{ij}-s_0b^ib^j-s_{-1}(y^ib^j+y^jb^i)-s_{-2}y^iy^j,
\end{eqnarray*}
\end{description}
where \vspace{-0.2cm}
 \begin{eqnarray*}
s_0:&=&e^{-\sigma}f^2q_0/(\varepsilon pL^2),
\hspace{.2cm}s_{-1}:=p_{-1}f^2/(p \,\varepsilon L^2),\hspace{.2cm}
s_{-2}:=p_{-1}(e^{\sigma}m^2 p L^2-b^2f^2)/(\varepsilon p \beta
L^2),
\end{eqnarray*}
   $\varepsilon :=f^2(e^\sigma p+m^2 q_0)/L^2\neq 0$,\quad
 $m^2=g^{ij}m_jm_i=m^im_i\neq 0$,\quad $b^i=g^{ij}b_j$.
  \end{prop}
   \begin{rem}\label{s01}
The quantities $s_{0}\,,\,s_{-1}\,,\,s_{-2}$
 satisfy:
$$\beta s_0+L^2s_{-1}=q/\varepsilon,$$
$$b^2s_{-1}+\beta s_{-2}=e^{\sigma}p_{-1}m^2/\varepsilon.$$
\end{rem}

\bigskip

Let $C_i=C_{ijk}g^{jk}$, $C^i=C^i_{jk}g^{jk}$ and $C^2=C^iC_i$.
Then, we have
\begin{prop}\label{c1} Under a generalized $\beta$-conformal change, we have
\begin{description}
    \item[(a)] The Cartan  tensor $\overline{C}_{ijk}$ has the form
    \begin{equation}\label{cc1}
    \overline{C}_{ijk}=e^{\sigma} p\,C_{ijk}+V_{ijk},
\end{equation}

    \item[(b)] The (h)hv-torsion tensor $\overline{C}^{l}_{ij}$ has the form
\begin{equation}\label{sh1-101}
   \overline{C}^l_{ij}=C^l_{ij}+M^l_{ij}\,,
\end{equation}
\end{description}
 where
\vspace{-.5cm}
\begin{eqnarray}\label{vijk}
    V_{ijk}&:=&\frac{e^{\sigma}p_{-1}}{2}(h_{ij}m_k+h_{jk}m_i+h_{ki}m_j)
    +\frac{p_{02}}{2}m_im_jm_k,\\
\nonumber \qquad\qquad M^l_{ij}&:=&\frac{1}{2}(e^{-\sigma}m^l/p-m^2
    (s_0b^l+s_{-1}y^l))(p_{02}m_im_j+e^{\sigma}p_{-1}h_{ij})\\
 \nonumber   &&-e^{\sigma}(s_0b^l+s_{-1}y^l)(p\,C_{ij\beta}+p_{-1}m_im_j)
 \nonumber  +\frac{p_{-1}}{2p}(h^l_im_j+h^l_jm_i),
\end{eqnarray}
and \, $h^i_j=g^{il}h_{lj}$, $C_{ij\beta}:=C_{ijr}b^r$,
$C_{i\beta\beta}:=C_{ijk}b^jb^k$ and so on.
\bigskip

\noindent\textbf{\emph{(c)}}  $\overline{C}_i=C_i-e^\sigma
p\,s_0C_{i\beta\beta}+\lambda m_i$,

\noindent where  \quad $
\displaystyle{\lambda:=\frac{(n+1)p_{-1}}{2p}-\frac{3e^\sigma
p_{-1}m^2 s_0}{2}+\frac{p_{02}m^2}{2(e^\sigma p+q_0 m^2)}},$
\bigskip

\noindent\textbf{\emph{(d)}} $\displaystyle{
\overline{C}^i=\frac{e^{-\sigma}}{p} \,C^i
    +J^i}$,  $$\text{where}\quad J^i:=\frac{\lambda e^{-\sigma}}{p}\, m^i-s_0C^i_{\beta\beta}-(C_\beta+\lambda m^2-
   e^{\sigma}s_0p\,C_{\beta\beta\beta})(s_0b^i+s_{-1}y^i), \,\, C_{\beta}:=C_ib^i.$$
    \bigskip

\noindent\textbf{\emph{(e)}}
$\displaystyle{\overline{C}^2=\frac{e^{-\sigma}}{p}\,C^2+\Phi}$,

\noindent where
    \vspace{-0.9cm}\begin{eqnarray*}
   \Phi &:=&  \lambda^2 m^2(( e^{-\sigma}/p)-s_0m^2)+C_{\beta}((2\lambda e^{-\sigma}/p)
    -s_0(1+2\lambda m^2))\phantom{then}\\
    &&+s_0C_{\beta\beta\beta}(1-3\lambda+e^{\sigma}s_0p \,C_\beta)\\
    &&
  +s_0 C_{r\beta\beta}(e^{2\sigma}s_0^2p^2C_{\beta\beta\beta}b^r-\lambda
  s_0
   m^2b^r-e^{\sigma}s_0p \,  C^r_{\beta\beta}-2C^r).
\end{eqnarray*}
\end{prop}
\begin{prop}\label{s}
 Under  a generalized     $\beta$-conformal change, the v-curvature
 tensor of $(M,L)$ is transformed as follows:
 \begin{eqnarray*}
    \overline{S}_{lijk}&=&e^\sigma p S_{lijk}
    +\mathfrak{A}_{jk}\{H_{lk}h_{ij}+H_{ij}h_{lk}+\omega_{lk}C_{ij\beta}+\omega_{ij}C_{lk\beta}\},
 \end{eqnarray*}

where
\begin{eqnarray*}
  H_{ij}: &=&K_1m_im_k+K_2C_{ij\beta}+K_3h_{ij},\quad
   \omega_{ij}:=K_4m_im_j-\frac{1}{2}e^{2\sigma}p^2s_0C_{ij\beta},\\
   K_1:&=&\frac{e^{2\sigma}p_{-1}^2}{4p}(e^{-\sigma}-2s_0 p\, m^2)
  +\frac{e^\sigma p_{-1}p_{02}m^2}{4(e^\sigma
p+q_0m^2)},\quad
  K_2:=\frac{e^\sigma p_{-1}}{2}-\frac{1}{2}e^{2\sigma}s_0p \,p_{-1}m^2,  \\
  K_3:&=&\frac{e^{2\sigma} p_{-1}^2m^2}{8 (e^\sigma p+q_0m^2)},\quad
  K_4:=\frac{e^\sigma p p_{02}}{2(e^\sigma p+q_0m^2)}-e^{2\sigma} s_0p\,
  p_{-1}.
\end{eqnarray*}
\end{prop}
\begin{rem}~\par The tensors $H_{ij}$ and $\omega_{ij}$ defined above
have the following properties:
\begin{description}
    \item[(1)] $H_{ij}$ and $\omega_{ij}$ are symmetric.

    \item[(2)] $H_{ij}$ and $\omega_{ij}$ are indicatory:    $H_{ij}y^i=0$, $\omega_{ij}y^i=0$.

    \item[(3)]$g^{ij}H_{ij}=K_1m^2+K_2C_\beta+(n-1)K_3$  and \,  $g^{ij}\omega_{ij}=K_4m^2
    -\frac{1}{2}e^{2\sigma}s_0p^2C_\beta$.
\end{description}
\end{rem}
 \begin{prop}\label{ss}
 Under a generalized     $\beta$-conformal change, the vertical Ricci tensor $\overline{S}_{ik}$
  and the vertical scalar curvature  $\overline{S}$ associated with the transformed  space $(M,\overline{L})$
  are given by:
 \begin{eqnarray*}
    \overline{S}_{ik}&=& S_{ik}
    +K h_{ik}+\left(s_0m^2-\frac{e^{-\sigma}}{(n-3)p}\right)H_{ik}+\Psi_{ik},\\
    \overline{S}&=&\frac{e^{-\sigma}}{p}S+\frac{2e^{-\sigma}}{p}K
    \{(n-2)-e^\sigma p s_0m^2\}
    -s_0\Psi_{\beta\beta}+\frac{e^{-\sigma}}{p}\Psi-s_0S_{ik}b^ib^k,
 \end{eqnarray*}

 where\begin{eqnarray*}
 K &:=&s_0H_{\beta\beta}-\frac{e^{-\sigma}}{ p}(K_1m^2+K_2C_\beta+(n-1)K_3),\\
\Psi_{ik}&:=&\frac{e^{-\sigma}}{p}\{\omega_{rk}C^r_{i\beta}+\omega_{ri}C^r_{k\beta}-(K_4m^2
-\frac{1}{2}e^{2\sigma}s_0p^2C_\beta)C_{ik\beta}\} -s_0\{H_{\beta
k}m_i+H_{i\beta}m_k\\&&+\omega_{\beta
k}C_{i\beta\beta}+\omega_{i\beta}C_{k\beta\beta}
-\omega_{\beta\beta}C_{ik\beta}-
\omega_{ik}C_{\beta\beta\beta}+e^\sigma p S_{hijk}b^hb^j\},\\
H_{\beta\beta} &:=&H_{ij}b^ib^j,\quad
\omega_{\beta\beta}:=\omega_{ij}b^ib^j, \quad
\Psi:=\Psi_{ij}g^{ij}, \quad \Psi_{\beta\beta}:=\Psi_{ij}b^ib^j.
\end{eqnarray*}
\end{prop}
Note that the tensor $\Psi_{ij}$ is symmetric and indicatory.

%%%%%%%%%%%%%%%%%%%%%%%%%%%%%%%%%%%%%%%%%%%%%% 2.Special Finsler Spaces %%%%%%%%%%%%%%%%%%%%%%%%%%%%%%%%%%%%%%%%%%%%%%%%%%%

\Section{\textbf{Special Finsler spaces}}

In this section we will investigate the effect of the generalized
$\beta$-conformal change (\ref{change}) on   some special Finsler
space. Some of the results  obtained in  this section are
generalizations of known results and some  are new. For a
systematic study of special Finsler spaces, we refer to
\cite{r2.25}.\\

In what follows, let $(M,L)$ be a Finsler  manifold and
$(M,\overline{L})$ the transformed Finsler manifold under a
generalized $\beta$-conformal change. The geometric objects
associated with $(M,\overline{L})$ will be denoted by barred
symbols.

\begin{thm}\label{energy} For
$n>2$, under  a generalized $\beta$-conformal change, the
following assertions are equivalent
\begin{description}
    \item[(a)]$p_{-1}=0.$
    \item[(b)]$q=k\beta$; $k$ is a nonzero constant.
    \item[(c)]$\overline{C}_{ijk}=e^\sigma p\, C_{ijk}$.
    \item[(d)]$\overline{L}=(k'e^{2\sigma}L^2+k\beta^2)^\frac{1}{2}$;
    $k'$ is a nonzero constant.
\end{description}
\noindent The special $\beta$-conformal change \emph{\textbf{(d)}}
is referred to
    as an energy $\beta$-change \cite{energy.}.
\end{thm}
\begin{proof}~\par
\noindent\textbf{(a)} $\Rightarrow$\textbf{(b)}:  Let  $p_{-1}=0,$
then $\displaystyle{\frac{f f_{12}}{L}+\frac{p f_2}{f}=0}$ which
leads to $ff_{12}+f_1f_2=0,$ hence,
$\displaystyle{\frac{\partial}{\partial \widetilde{L}}(ff_2)=0}$. By
integration, taking the homogeneity of $f$ into account, we get $q=k
\beta$, with  $k\neq 0$.\\

\noindent\textbf{(b)}$\Rightarrow$\textbf{(c)}:  Let $q=k\beta$,
then $\displaystyle{\frac{\partial}{\partial
\widetilde{L}}(ff_2)=ff_{12}+f_1f_2=0}$, which leads to $p_{-1}=0$.
Using  $\beta p_0+e^\sigma L^2p_{-1}=q$, we get  $\beta p_0=q$. By
differentiating the last identity with respect to $\beta$, we have
$$\beta p_{02}+p_0=f_2^2+ff_{22}=p_0,$$ which leads to $p_{02}=0$.
Hence, by (\ref{vijk}) $V_{ijk}=0$ and, consequently,
$\overline{C}_{ijk}=e^\sigma p
C_{ijk}$.\\

\noindent\textbf{(c)}$\Rightarrow$\textbf{(d)}:  Let $V_{ijk}=0$,
then
$$e^{\sigma}p_{-1}(h_{ij}m_k+h_{jk}m_i+h_{ki}m_j)+p_{02}m_im_jm_k=0.$$
By  contraction by $b^i$, we have
\begin{equation}\label{proof}
e^{\sigma}p_{-1}(2m_jm_k+m^2h_{jk})+p_{02}m^2m_jm_k=0.
\end{equation}
 Contracting  (\ref{proof}) again  by $b^j$, we get $3 e^\sigma p_{-1}=m^2
 p_{02}$.
Hence, (\ref{proof}) reduces to $p_{-1}(m^2h_{jk}-m_jm_k)=0,$ which
leads to $p_{-1}=0$ or $m^2h_{jk}-m_jm_k=0$. Now, if
$m^2h_{jk}-m_jm_k=0$, then, $n=2$ which contradicts the hypothesis.
Hence, $p_{-1}=0$, and consequently,  $q=k\beta$. Then, we have the
partial differential equation
$$ff_2=k\beta.$$
By integration with respect to $\beta$ and using the fact that $f$
is homogenous of degree 1 in $\beta$ and $ \widetilde{L}$, we get
$$f^2=k\beta^2+\varphi( \widetilde{L}),$$
where $\varphi(\widetilde{L})$ is a homogenous function of degree 2
in $\widetilde{L}$, which may be written as
$\varphi(\widetilde{L})=k'\widetilde{L}^2$. Hence,
$f^2=k\beta^2+k'\widetilde{L}^2$ and consequently,
$$\overline{L}=(k'\widetilde{L}^2+k\beta^2)^\frac{1}{2}
=(k'e^{2\sigma}L^2+k\beta^2)^\frac{1}{2}.$$

 \noindent\textbf{(d)}$\Rightarrow$\textbf{(a)}:  It is obvious.
\end{proof}

\begin{cor} For $n>2$,  under  a generalized $\beta$-conformal
change, if one of the above equivalent conditions holds, then the
space $(M,\overline{L})$ is  Riemannian if and
    only if $(M,L)$ is  Riemannian.
\end{cor}

 We will study the change of some special Finsler spaces under a
 generalized $\beta$-conformal
change.

\begin{defn}
A Finsler space $(M,L)$ with dimension  $n\geq 3$   is said to be
quasi-C-reducible   if the Cartan   tensor $C_{ijk}$ satisfies
\begin{equation}\label{paper2.8}
C_{ijk}=Q_{ij}C_k+Q_{jk}C_i+Q_{ki}C_j,
\end{equation}
where $Q_{ij}$ is a symmetric indicatory  tensor.
\end{defn}
By    Proposition  \ref{c1}, assuming $\lambda \neq 0$, we have
\begin{eqnarray*}
  \overline{C}_{ijk} &=&e^{\sigma} p\,C_{ijk}+\frac{e^{\sigma}p_{-1}}{2}(h_{ij}m_k+h_{jk}m_i+h_{ki}m_j)
    +\frac{p_{02}}{2}m_im_jm_k  \\
   &=& e^{\sigma} p\,C_{ijk}+\frac{1}{6}\mathfrak{S}_{ijk}\{(3e^{\sigma}p_{-1}h_{ij}+p_{02}m_im_j)m_k\} \\
   &=&  e^{\sigma} p\,C_{ijk}+\frac{1}{6 \lambda}\mathfrak{S}_{ijk}\set{(3e^{\sigma}p_{-1}h_{ij}
   +p_{02}m_im_j)(\overline{C}_k-C_k+e^\sigma p s_0 C_{k\beta\beta})} \\
   &=& e^{\sigma} p\,C_{ijk}+ \frac{1}{6
   \lambda}\mathfrak{S}_{ijk}\set{(3e^{\sigma}p_{-1}h_{ij}+p_{02}m_im_j)\overline{C}_k}\\
   && + \frac{1}{6
   \lambda}\mathfrak{S}_{ijk}\set{(3e^{\sigma}p_{-1}h_{ij}+p_{02}m_im_j)(e^\sigma p
   s_0C_{k\beta\beta}-C_k)}.
\end{eqnarray*}
Hence, we have
\begin{lem}\label{quasi1} Under a generalized     $\beta$-conformal change, the transformed
 Cartan tensor can be written in the form

 $$\overline{C}_{ijk}=
 \mathfrak{S}_{ijk}\set{\overline{Q}_{ij}\overline{C}_k}+q_{ijk},$$

where $\overline{Q}_{ij}:=\frac{1}{6
\lambda}(3e^{\sigma}p_{-1}h_{ij}
   +p_{02}m_im_j),$
   $$q_{ijk}:=
   \frac{1}{6 \lambda}\mathfrak{S}_{ijk}\set{2e^\sigma \lambda p\, C_{ijk} +(3e^{\sigma}p_{-1}h_{ij}+p_{02}m_im_j)(e^\sigma p
   s_0C_{k\beta\beta}-C_k)}.$$

   \end{lem}
  By the above lemma and  taking into account that  the tensor $\overline{Q}_{ij}$ is symmetric and indicatory,
       we get the following result.
\begin{thm}\label{paper2.9}  If the tensor $q_{ijk}=0$, then the
space $(M, \overline{L})$ is  quasi-C-reducible.
\end{thm}

As a corollary  of the above theorem, we obtain a generalized form
of Matsumoto's result \cite{r2.22}:

\begin{cor}\label{paper2.10}Under a generalized     $\beta$-conformal change,
 a Reimannian space $(M,L)$ is
transformed to a quasi-C-reducible space.
\end{cor}

\begin{defn}\label{semi-1}
A Finsler space $(M,L)$ of dimension $n\geq 3$ is called
semi-C-reducible,  if the Cartan  tensor $C_{ijk}$ is written in the
form:
\begin{equation}\label{semi-c-red.}
C_{ijk}=\frac{r}{n+1}(h_{ij}C_k+h_{ki}C_j+h_{jk}C_i)+\frac{t}{C^2}C_iC_jC_k,
\end{equation}
 where $r$ and $t$ are scalar functions such that   $r+t=1$.
\end{defn}
The next result has been  obtained by Matsumoto and Shibata
\cite{r2.21} in  the special case of Finsler spaces with $(\alpha,
\beta)$-metric.
\begin{thm}\label{semi-4} A Riemannian space is transformed to a
semi-C-reducible space,
   by a generalized
$\beta$-conformal change.
\end{thm}
\begin{proof} From Proposition \ref{h-g.1} and  Proposition \ref{c1}, we get
\begin{eqnarray*}
   \overline{C}_{ijk} &=&\frac{1}{2}e^{\sigma}p_{-1}(h_{ij}m_k+h_{jk}m_i+h_{ki}m_j)+\frac{1}{2}p_{02}m_im_jm_k \\
   &=&\frac{p_{-1}}{2p\lambda}(\overline{h}_{ij}\overline{C}_k
   +\overline{h}_{jk}\overline{C}_i
  + \overline{h}_{ki}\overline{C}_kj)+\frac{m^2(pp_{02}-3p_{-1}q_0)}{2p\lambda(e^\sigma p
  +m^2q_0)\overline{C}^2}\overline{C}_i\overline{C}_j\overline{C}_k  \\
   &=&\frac{r}{n+1}(\overline{h}_{ij}\overline{C}_k
   +\overline{h}_{jk}\overline{C}_i
  +
  \overline{h}_{ki}\overline{C}_kj)+\frac{t}{\overline{C}^2}\overline{C}_i\overline{C}_j\overline{C}_k,
\end{eqnarray*}
where $$r=\frac{p_{-1}(n+1)}{2p\lambda},\quad
t=\frac{m^2(pp_{02}-3p_{-1}q_0)}{2p\lambda(e^\sigma p+m^2q_0)},
\quad r+t=1,$$ which means that $(M,\overline{L})$ is
semi-reducible.
\end{proof}

\begin{defn}
A Finsler space $(M,L)$ of dimension $n\geq 3$ is called C-reducible
if the Cartan tensor $C_{ijk}$ has the form:
\begin{equation}\label{c-reducible*}
C_{ijk}=h_{ij}A_k+h_{ki}A_j+h_{jk}A_i, \quad A_i=\frac{C_i}{n+1}.
\end{equation}
\end{defn}

Define the  tensor
$$K_{ijk}=C_{ijk}-(h_{ij}A_k+h_{ki}A_j+h_{jk}A_i).$$

It is clear that $K_{ijk}$ is symmetric and indicatory. Moreover,
$K_{ijk}$ vanishes if and only if the Finsler space $(M,L)$ is
C-reducible.

\begin{prop}\label{c-1}
Under a generalized     $\beta$-conformal change, the tensor
$\overline{K}_{ijk}$ associated with  the space $(M,\overline{L})$
has the form
\begin{eqnarray*}
   \overline{K}_{ijk}&=& e^\sigma p
   K_{ijk}+d_{ijk},
\end{eqnarray*}
where
\begin{eqnarray*}
d_{ijk}&:=&\frac{1}{n+1}\mathfrak{S}_{ijk}\{(n+1)(\alpha_1 h_{ij}+
   \alpha_2 m_im_j)m_k+q_0 m_im_jC_k\\
   &&+(s_0 p q_0 m_im_j+e^\sigma p^2 s_0
   h_{ij})C_{k\beta\beta}\},\\
 \alpha_1&:=&\frac{e^\sigma p_{-1}}{2}-\frac{e^\sigma p\lambda}{n+1},\quad
   \alpha_2:=\frac{ p_{02}}{6}-\frac{q_0\lambda}{n+1}.
\end{eqnarray*}
\end{prop}

Consequently, we have
\begin{thm}\label{c-reducible}   Under a generalized     $\beta$-conformal change, the
following assertions
\begin{description}
    \item[(a)]the space $(M,L)$ is C-reducible,
    \item[(b)]the space $(M,\overline{L})$ is C-reducible
\end{description}
are equivalent if and only if the tensor $d_{ijk}$ vanishes.
\end{thm}

\begin{cor}If
$\overline{L}=e^\sigma L+\beta$, $L$ being Finslerian,  then the
tensor $d_{ijk}$ vanishes. Consequently,  $(M,\overline{L})$ is
C-reducible if and only if  $(M,L)$ is C-reducible.
\end{cor}

\begin{lem}\label{lr-1} Under
a generalized $\beta$-conformal change $\overline{L}=f(e^\sigma
L,\beta)$, with $L$ Riemannian, the tensor $d_{ijk}$ takes the form
$$d_{ijk}=\mathfrak{S}_{ijk}\{\alpha_1 h_{ij}m_k+
   \alpha_2 m_im_jm_k\}.$$
\end{lem}

\begin{thm}\label{eqv.}
Under  a generalized $\beta$-conformal change
$\overline{L}=f(e^\sigma L,\beta)$, with $L$ Riemannian, the
following assertions are equivalent:
\begin{description}
    \item[(a)]$\alpha_1=0$ and $\alpha_2=0$,
    \item[(b)]$(M,\overline{L})$ is C-reducible,
    \item[(c)]$(M,\overline{L})$ is either  of Randers type or of Kropina
    type.
\end{description}
\end{thm}
\begin{proof}~\par
\noindent\textbf{(a)} $\Rightarrow$\textbf{(b)}:  It is obvious.

 \noindent\textbf{(b)} $\Rightarrow$\textbf{(a)}:  Let the space
 $(M,\overline{L})$ be  C-reducible, then, by Lemma
\ref{lr-1}, we have:
  \begin{equation}\label{proof2}
d_{ijk}= \mathfrak{S}_{ijk}\{\alpha_1 h_{ij}m_k+
   \alpha_2 m_im_jm_k\}=0.
\end{equation} Contracting (\ref{proof2}) by $g^{ij}$, we get
\begin{equation}\label{proof3}
(n+1)\alpha_1+3m^2\alpha_2=0,
\end{equation} and  contracting the same equation by  $b^ib^j$, we
get  \begin{equation}\label{proof3} \alpha_1+m^2\alpha_2=0.
\end{equation} The last
two relations lead to $(n-2)\alpha_1=0$. Since $n>2$, then
$\alpha_1=0$ and consequently $\alpha_2=0$, by (\ref{proof3}).

\noindent\textbf{(a)} $\Rightarrow$\textbf{(c)}:   If
$\alpha_1=\alpha_2=0$, we have $(n+1)p_{-1}=2p\lambda$ and
$(n+1)p_{02}=6q_0\lambda$. Solving the last two equations for
$\lambda$, we get
$$3 q_0 p_{-1}=pp_{02}.$$ From which
we obtain the  patrial differential equation
$$\frac{3 f_{12}f_{22}}{f_1}-f_{222}=0.$$ Now, if $f_{22}=0$, by
integration with respect to $\beta$ and taking  the homogeneity
  of $f$ into account, we get $f_2=\varphi_1(
\widetilde{L})$, where $\varphi_1( \widetilde{L})$ is a homogenous
function of degree $0$ in $\widetilde{L}$. Hence, by integrating
$f_2$ with respect to $\beta$, we get
$$\overline{L}=\varphi_1( \widetilde{L}) \beta+\varphi_2( \widetilde{L}),$$
where $\varphi_2( \widetilde{L})$ is a homogenous function of degree
$1$ in $\widetilde{L}$. By the homogeneity properties of $\varphi_1(
\widetilde{L})$ and $\varphi_2( \widetilde{L})$, using Euler
theorem, we conclude that $\varphi_1( \widetilde{L})=c_1$ and
$\varphi_2( \widetilde{L})=c_2$, where $c_1$ and $c_2$ are
constants. Consequently,
$$\overline{L}=c_2
\widetilde{L} +c_1 \beta.$$

On the other hand, if $f_{22}\neq 0$,  we have $$\frac{3
f_{12}}{f_1}-\frac{f_{222}}{f_{22}}=0,$$
 which, by
integration with respect to $\beta$, gives

$$\displaystyle{3\ln f_1-\ln f_{22}=\ln \varphi_3( \widetilde{L})}\Rightarrow
\displaystyle{\frac{f_1^3}{f_{22}}=\varphi_3( \widetilde{L})=c_3
 \widetilde{L}}, $$
where $\varphi_3( \widetilde{L})$ is a homogenous   function of
degree 1 in  $\widetilde{L}$ and $c_3$ is nonzero constant.  Using
$\widetilde{L} f_{11} +\beta f_{12}=0$ and $\widetilde{L} f_{21}
+\beta f_{22}=0$, we have
$$\frac{f_{11}}{f_1^3}=\frac{c_3 \beta^2}{\widetilde{L}^3},$$ from
which $\displaystyle{f_1=\frac{\widetilde{L}}{\sqrt{c_3 \beta^2+ c_4
\widetilde{L}^2}}}$. If $c_4\neq 0$, then
$$f=\frac{1}{c_4}\sqrt{c_3 \beta^2+ c_4 \widetilde{L}^2}+c_5 \beta
,$$ and if  $c_4=0 $, then
$$f=\frac{\widetilde{L}^2+\beta^2}{c_3 \beta}.$$
The former may be regarded as of  Randers type and the later as of
Kropina type.

\noindent\textbf{(c)} $\Rightarrow$\textbf{(a)}:  The result follows
directly by computing $\alpha_1$ and $\alpha_2$ for Randers and
Kropina spaces.
\end{proof}

It should be noted  that Matsumoto \cite{r2.9} showed that
C-reducible Finsler spaces with $(\alpha,\beta)$-metric are either
of Randers type or of Kropina type.

\begin{defn}
A Finsler space $(M,L)$ of dimension $n\geq 2$ is said to be
$C_2$-like  if the Cartan  tensor $C_{ijk}$ satisfies
\begin{equation}\label{c2-likedef.}
C^2C_{ijk}=C_iC_jC_k.
\end{equation}
\end{defn}

Let us  define the  tensor
$$\eta_{ijk}=C^2C_{ijk}-C_iC_jC_k.$$

It is clear that $\eta_{ijk}$ is symmetric and indicatory.
Moreover, $\eta_{ijk}$ vanishes if and only if the Finsler space
is $C_2$-like space.

\begin{prop}\label{c-1}
Under a generalized     $\beta$-conformal change, the tensor
$\overline{\eta}_{ijk}$ associated with  the space
$(M,\overline{L})$  has the form
\begin{eqnarray*}
   \overline{\eta}_{ijk}&=&
   \eta_{ijk}+I_{ijk},
\end{eqnarray*}
where \begin{eqnarray*}
 I_{ijk}&:=&(e^{-\sigma}/p)C^2 V_{ijk}+\Phi (e^\sigma p C_{ijk}+V_{ijk})-\lambda^3 m_im_jm_k
 \\
 &&-\lambda^2( m_jm_kC_i+m_im_jC_k+m_im_kC_j)-\lambda( m_kC_iC_j+m_jC_kC_i+m_iC_jC_k)\\
 &&- e^\sigma ps_0[C_{i\beta\beta}\{ e^\sigma p s_0C_{k\beta\beta}( \lambda m_j
 -C_j-e^{\sigma} p s_0 C_{j\beta\beta}) + e^\sigma
  p s_0C_{j\beta\beta}( \lambda m_k
 +C_k)\\
 &&-\lambda (m_kC_j
 + m_jC_k-\lambda m_jm_k)-C_jC_k\}
   +C_{k\beta\beta}(C_iC_j
-e^\sigma ps_0C_iC_{j\beta\beta}
 +\lambda m_iC_j\\&&+\lambda m_jC_i
 +\lambda^2 m_im_j)+\lambda C_{j\beta\beta}(\lambda m_im_k+C_iC_k
  + m_iC_k+ m_kC_i-e^\sigma p s_0 m_iC_{k\beta\beta})].
\end{eqnarray*}
\end{prop}

\begin{thm}\label{c2-like}   Under
a generalized     $\beta$-conformal change, the following
assertions
\begin{description}
    \item[(a)]the space $(M,L)$ is $C_2$-like,
    \item[(b)]the space $(M,\overline{L})$ is $C_2$-like
\end{description}
are equivalent if and only if the tensor $I_{ijk}$ vanishes.
\end{thm}
\begin{lem}\label{c2-like2}Starting with a Riemannian space $(M,L)$,  under
a generalized  $\beta$-conformal change, the tensor $I_{ijk}$ takes
the form:
\begin{equation}\label{sense}
I_{ijk}=\Phi V_{ijk}-\lambda^3m_im_jm_k.
\end{equation}

\end{lem}
\begin{thm}\label{c2-likenew}  For a $\beta$-conformal change
 $\overline{L}=e^\sigma L +\beta$; L being  Finslerian, a
necessary condition for the assertions
\begin{description}
    \item[(a)]the space $(M,L)$ is $C_2$-like,
    \item[(b)]the space $(M,\overline{L})$ is $C_2$-like
\end{description}
to be equivalent is that $C_\beta=0$.
\end{thm}

\begin{proof}
In the case of $\overline{L}=e^\sigma L +\beta$; L being Finslerian,
$\Phi= \lambda e^{-\sigma}L(\lambda m^2+2C_\beta)/\overline{L}$,\\
 $\lambda= \frac{n+1}{2\overline{L}}$ and  $
V_{ijk}=\frac{e^{\sigma}}{2L}(h_{ij}m_k+h_{jk}m_i+h_{ki}m_j)$. Now,
let the above assertions be equivalent,  so $I_{ijk}=0$. Contracting
(\ref{sense}) by $g^{jk}$, we have
$\frac{n+1}{2\overline{L}}C_\beta=0$ and  the result follows.
\end{proof}

  If $(M,L)$ is a Riemannian space and the tensor $I_{ijk}$
  vanishes, i.e., $(M,\overline{L})$ is $C_2$-like, we have
  $$\lambda^2\frac{m^2}{e^\sigma p+q_0m^2}(e^{\sigma}p_{-1}(h_{ij}m_k+h_{jk}m_i+h_{ki}m_j)
    +p_{02}m_im_jm_k)-2\lambda^3m_im_jm_k=0,$$
    contracting by $b^ib^j$ and assuming that $\lambda\neq 0$, we
    get
    \begin{equation}\label{n2}
(n-2)p_{-1}=0.
\end{equation}
Hence, we have
\begin{thm}\label{paper2.6}
Starting with a Riemannian space $(M,L)$, if the transformed space
$(M,\overline{L})$ is $C_2$-like, then one of the following holds:
\begin{description}
    \item[(a)] $dim\, M=2$.
    \item[(b)] The generalized     $\beta$-conformal change is an energy $\beta$-change
    and  the transformed space is Rimannian.
\end{description}

\end{thm}
\begin{cor}\label{c2-like4}Let the generalized     $\beta$-conformal change
 be of the form $\overline{L}=e^\sigma
 L+\beta$, with L  Riemannian. If $(M,\overline{L})$ is $C_2$-like, then $dim\, M=2$.
\end{cor}
\begin{cor}\label{paper2.7}  A Reimannian space of dimension  $\geq 3$
can not be transformed to a non-Reimannian $C_2$-like space.
\end{cor}
Now, we are going to  study two special Finsler spaces whose
defining property  depends on the v-curvature tensor $S_{lijk}$,
namely, the $S_3$-like and $S_4$-like Finsler spaces.

\begin{defn}
A Finsler space $(M^n,L)$ with  dimension $n > 3$   is said to be
$S_3$-like  if the v-curvature  tensor $S_{lijk}$ satisfies
\begin{equation}\label{s3-likedef.}
S_{lijk}=\frac{S}{(n-1)(n-2)}\{h_{ik}h_{lj}-h_{ij}h_{lk}\},
\end{equation}
 where S  is the vertical scalar curvature.

\end{defn}

Define the following tensor
$$\mu_{lijk}=S_{lijk}-\frac{S}{(n-1)(n-2)}\{h_{ik}h_{lj}-h_{ij}h_{lk}\}.$$
It is clear that the tensor $\mu_{hijk}$ vanishes if and only if
the space is $S_3$-like.
\begin{prop}\label{c-1}
Under a generalized     $\beta$-conformal change, the tensor
$\overline{\mu}_{hijk}$ associated with  the space
$(M,\overline{L})$  has the form:
$$\overline{\mu}_{lijk} =e^\sigma p \mu_{lijk}+r_{lijk},$$
where
\begin{eqnarray*}
  r_{lijk} &=&
  \mathfrak{A}_{jk}\{H_{lk}h_{ij}+H_{ij}h_{lk}+\omega_{lk}C_{\beta ij}+\omega_{ij}C_{\beta lk}
  -\frac{e^{2\sigma}p^2\Omega}{(n-1)(n-2)}h_{lj}h_{ik}\\
   &&- \frac{q_0}{(n-1)(n-2)}(S+e^{\sigma}p\Omega)
    (h_{ik}m_lm_j+h_{lj}m_im_k)\},\\
   \Omega
   &:=&\frac{e^{-\sigma}}{p}\Psi-s_0S_{ik}b^ib^k-s_0\Psi_{\beta\beta}
   +\frac{2e^{-\sigma}}{p}K(n-2-e^\sigma p s_0m^2).
\end{eqnarray*}
\end{prop}
\begin{thm}\label{s3like}   Under a generalized     $\beta$-conformal change, the
following assertions
\begin{description}
    \item[(a)]the space $(M,L)$ is $S_{3}$-like,
    \item[(b)]the space $(M,\overline{L})$ is $S_{3}$-like
\end{description}
are equivalent if and only if the tensor $r_{lijk}$ vanishes.
\end{thm}
\begin{prop}\label{sabd} For a $\beta$-conformal change $\overline{L}=e^\sigma
L+\beta$,  the tensor $r_{lijk}$ takes the form
$$r_{lijk}=C_{\beta jk}h_{il}+\frac{1}{2\overline{L}}m_jm_kh_{il}+\frac{m^2}{4\overline{L}}h_{jk}h_{il}
-A_\beta h_{jk}h_{il},$$ where
$A_\beta=\frac{1}{n-1}\left(C_\beta+\frac{n+1}{4\overline{L}}m^2\right)$.
\end{prop}
From the above proposition, we retrieve a result due to  Abed
\cite{r2.1}:
\begin{cor}
In the case of a  $\beta$-conformal change $\overline{L}=e^\sigma
L+\beta$, the following assertions
\begin{description}
    \item[(a)]the space $(M,L)$ is $S_{3}$-like,
    \item[(b)]the space $(M,\overline{L})$ is $S_{3}$-like
\end{description}
are equivalent if and only if
$$C_{rjk}b^r+\frac{1}{2\overline{L}}m_jm_k+\frac{m^2}{4\overline{L}}h_{jk}
=A_\beta h_{jk}.$$
\end{cor}
\bigskip
 Finally, we study  $S_4$-like Finsler spaces under a
generalized $\beta$-conformal change.

\begin{defn}
A Finsler space $(M,L)$ with  dimension $n > 4$   is said to be
$S_4$-like  if the v-curvature  tensor $S_{hijk}$ satisfies
\begin{equation}\label{paper2.16}
S_{lijk}=\mathfrak{A}_{jk}\{h_{lj}\mathbb{M}_{ik}+h_{ik}\mathbb{M}_{lj}\},
\end{equation}
where
$\displaystyle{\mathbb{M}_{ik}=\frac{1}{(n-3)}\set{S_{ik}-\frac{Sh_{ik}}{2(n-2)}}}.$
\end{defn}

Define the  tensor
$$\zeta_{lijk}=S_{lijk}-\mathfrak{A}_{jk}\{h_{lj}\mathbb{M}_{ik}+h_{ik}\mathbb{M}_{lj}\}.$$
It is clear that the tensor $\zeta_{hijk}$ vanishes if and only if
the  space is $S_4$-like.
\begin{prop}\label{s4}
Under a generalized generalized     $\beta$-conformal change, the
tensor $\overline{\zeta}_{hijk}$  associated with  the space
$(M,\overline{L})$ has the form
 $$ \overline{\zeta}_{lijk} =e^\sigma p
  \zeta_{lijk}+\varepsilon_{lijk},$$
  where
\begin{eqnarray*}
 \varepsilon_{lijk}:&=&\mathfrak{F}_{lijk}\{\omega_{lk}C_{ijs}b^s+q_0
 \mathbb{M}_{lk}m_im_j
  +\frac{e^{\sigma}p}{n-3}(s_0m^2H_{lk}h_{ij}+Kh_{lk}h_{ij}
   +\Psi_{lk}h_{ij}) \\&&-\frac{e^{\sigma}}{n-3}\left( \frac{e^{-\sigma} q_0 S
   }{2(n-2)}m_lm_kh_{ij} +\frac{\Omega p (e^\sigma
   ph_{lk}+q_0m_lm_k)}{2(n-2)}h_{ij}+e^{-\sigma }q_0Kh_{lk}m_im_j\right)\\&&+\frac{q_0}{(n-3)p}\left
   (\Psi_{lk}+(s_0 p m^2-(n-3)e^{-\sigma})H_{lk}
   -\frac{\Omega e^\sigma p^2}{2(n-2)}h_{lk}\right)m_im_j\},
\end{eqnarray*}
and
$$\mathfrak{F}_{lijk}\{X_{lk}Y_{ij}\}:=X_{lk}Y_{ij}+X_{ij}Y_{lk}-X_{lj}Y_{ik}-X_{ik}Y_{lj}.$$
\end{prop}
\begin{thm}\label{s4like}   Under a generalized     $\beta$-conformal change, the
following assertions
\begin{description}
    \item[(a)]the space $(M,L)$ is $S_{4}$-like,
    \item[(b)]the space $(M,\overline{L})$ is $S_{4}$-like
\end{description}
are equivalent if and only if the tensor $\varepsilon_{hijk}$
vanishes.
\end{thm}
In the case of a $\beta$-conformal change $\overline{L}=e^\sigma
L+\beta$, the tensor $\varepsilon_{hijk}$ vanishes and we retrieve
the the following result of Abed \cite{r2.1}.
\begin{cor}\label{sabd2} For a  of $\beta$-conformal change $\overline{L}=e^\sigma
L+\beta$,  the space $(M,L)$ is $S_4$-like
 if and only if  the space
$(M,\overline{L})$ is $S_{4}$-like.\\
\end{cor}

%%%%%%%%%%%%%%%%%%%%%%%%%%%%%%%%%%%%%%%%%%%%%%%%%%%% 3. The T-tensor %%%%%%%%%%%%%%%%%%%%%%%%%%%%%%%%%%%%%%%%%%%%%%%%%%%

\Section{\textbf{The  T-tensor  $T_{hijk}$}}
The T-tensor is defined by \cite{ttensor}
$$T_{hijk}=LC_{hij}{\mid}_k
+C_{hij}l_k+C_{hik}l_j +C_{hjk}l_i+C_{ijk}l_h,$$
It should be noted that the T-tensor has a great contribution  in
geometric properties of special Finsler spaces.  For instance,
Hashiguchi \cite{r2.6} has shown that a Landsberg space remains
Landsberg under a conformal transformation, if and only if $
T_{hijk} = 0$. On the other hand, Matsumoto \cite{matsumoto} has
obtained interesting results for spaces with $T_{hijk} = 0$ and,
further, he investigated the
three-dimensional Finsler spaces with vanishing T-tensor.

In this section we compute the T-tensor under a generalized $\beta$-conformal change and consider some interesting special cases.

\begin{thm}\label{t}
Under a generalized $\beta$-conformal change, the transformed
T-tensor takes the form:
\begin{eqnarray*}
  \overline{T}_{lijk} &=& \frac{e^\sigma p \overline{L}}{L} T_{lijk} -\overline{L}
  (\frac{\beta e^\sigma
  p_{-1}}{2L^2}+2K_3)(h_{li}h_{jk}+h_{lj}h_{ik}+h_{lk}h_{ij})\\
  &&+(h_{li}\nu_{jk}+h_{lj}\nu_{ik}+h_{ij}\nu_{lk}
   +h_{jk}\nu_{il}+h_{lk}\nu_{ij}+h_{ik}\nu_{jl})\\
   &&+(e^\sigma p f_2-\frac{1}{2}e^\sigma  \overline{L} p_{-1})(C_{lij}m_k+
   C_{ijk}m_l+C_{jlk}m_i+C_{lik}m_j) \\
   &&-\overline{L}(M_{ij}C_{lk\beta}+M_{jl}C_{ik\beta}+M_{il}C_{jk\beta}
   +M_{lk}C_{ij\beta}+M_{jk}C_{il\beta}+M_{ik}C_{jl\beta})\\
   &&+\overline{L}e^{2\sigma}s_0p^2(C_{ij\beta}C_{lk\beta}+C_{lj\beta}C_{ik\beta}+C_{il\beta}C_{jk\beta})
   +\frac{1}{2}\overline{L}(6K_5+p_{022})m_lm_im_jm_k\\
   &&-\frac{\overline{L}p_{02}}{2L}(n_{ij}m_km_l
  +n_{lk}m_jm_j)
  +\frac{1}{2}p_{02}(\dot{n}_{ij}m_km_l
  +\dot{n}_{lk}m_im_j)
\end{eqnarray*}
where
\begin{eqnarray*} \nu_{ij}&:=&\frac{1}{2}e^\sigma
p_{-1}\dot{n}_{ij}-\overline{L}(K_1+\frac{3e^\sigma
p^2_{-1}}{4p}+\frac{\beta
p_{02}}{2L^2})m_im_j-\frac{\overline{L}e^\sigma p_{-1}}{2L}n_{ij},\\
K_5&:=&e^{2\sigma}s_0 p_{-1}^2-\frac{4 e^{\sigma}
p_{-1}p_{02}+p^2_{02}m^2}{ 4(e^\sigma p+q_0m^2)},\\
\dot{n}_{ij}&:=&\overline{l}_im_j+\overline{l}_jm_i.
\end{eqnarray*}
\end{thm}
\begin{proof}
One can show that
\begin{eqnarray}\label{lamda}
  \dot{\partial}_k\overline{C}_{lij}  &=& e^\sigma p \dot{\partial}_k C_{lij}
 \nonumber   +e^\sigma p_{-1}(C_{lij}m_k+C_{ijk}m_l+C_{jlk}m_i+C_{lik}m_j)\\
\nonumber    && -\frac{e^\sigma
p_{-1}}{2L}(h_{li}n_{jk}+h_{lj}n_{ik}+h_{ij}n_{lk}
   +h_{jk}n_{li}+h_{lk}n_{ij}+h_{ik}n_{jl}) \\
\nonumber    &&-\frac{\beta e^\sigma
p_{-1}}{2L^2}(h_{li}h_{jk}+h_{lj}h_{ik}+h_{lk}h_{ij})
   -\frac{\beta p_{02}}{2L^2}(h_{ij}m_km_l+h_{li}m_jm_k+h_{kl}m_im_j\\
\nonumber    &&+h_{jk}m_lm_i+h_{ik}m_jm_l+h_{lj}m_im_k)+\frac{1}{2}p_{022}m_im_jm_km_l\\
   &&-\frac{p_{02}}{2L}(n_{ij}m_km_l+n_{lk}m_im_j)
\end{eqnarray}
where  $n_{ij}:=l_im_j+l_jm_i$ and $p_{022}:=\frac{\partial
}{\partial\beta}p_{02}$. Similarly,
\begin{eqnarray}\label{cc}
\nonumber  \overline{C}_{ijr} \overline{C}^r_{lk} &=& e^\sigma p
  C_{ijr}C^r_{lk}
  +\frac{1}{2}e^\sigma p_{-1}(C_{ljk}m_i+C_{ilk}m_j+C_{ijk}m_l+C_{ijl}m_k)\\
\nonumber &&+K_4(C_{lk\beta}m_im_j+C_{ij\beta}m_lm_k)+
   (K_1+\frac{1}{4p}e^\sigma p^2_{-1})
   (h_{ij}m_km_l+h_{lk}m_im_j)\\
 \nonumber  &&+2K_3h_{ij}h_{lk}
   +K_2(C_{lk\beta}h_{ij}+C_{ij\beta}h_{lk})-K_5m_im_jm_lm_k-e^{2\sigma} p^2 s_0C_{ij\beta}C_{lk\beta}\\
   &&+\frac{e^\sigma
   p^2_{-1}}{4p}(h_{li}m_jm_k+h_{jk}m_hm_i+h_{ik}m_jm_h+h_{jl}m_im_k),
\end{eqnarray}

Using (\ref{lamda}) and (\ref{cc}),  we get
\begin{eqnarray}\label{ttensor2}
  \nonumber \overline{C}_{lij}\overline{|}_k &=&\dot{\partial}_k\overline{C}_{lij}- \overline{C}^m_{lk}\,\overline{C}_{mij}
  -\overline{C}^m_{ik}\,\overline{C}_{mlj}-\overline{C}^m_{jk}\,\overline{C}_{mli}\\
  \nonumber   &=&e^\sigma p\; C_{lij}|_k-(\frac{\beta e^\sigma
  p_{-1}}{2L^2}+2K_3)(h_{li}h_{jk}+h_{lj}h_{ik}+h_{lk}h_{ij})+3K_5m_im_jm_km_l
  \\
  \nonumber &&+\frac{1}{2}e^\sigma
  p_{-1}(C_{ljk}m_i+C_{ilk}m_j+C_{ijk}m_l+C_{ijl}m_k)-\frac{p_{02}}{2L}(n_{ij}m_km_l
  +n_{lk}m_jm_j)\\
  \nonumber &&
   -\frac{e^\sigma p_{-1}}{2L}(h_{li}n_{jk}+h_{lj}n_{ik}+h_{ij}n_{lk}
   +h_{jk}n_{lk}+h_{ij}n_{lk}+h_{ik}n_{jl})  \\
 \nonumber   &&-(K_1+\frac{3e^\sigma p^2_{-1}}{4p}+\frac{\beta p_{02}}{2L^2})(h_{ij}m_lm_k+h_{lk}m_im_j+h_{li}m_jm_k
   +h_{jk}m_lm_i\\
  \nonumber  &&+h_{ik}m_jm_l+h_{jl}m_im_k)-(M_{ij}C_{lk\beta}+M_{jl}C_{lik\beta}+M_{il}C_{jk\beta}
   +M_{lk}C_{ij\beta}\\&&+M_{jk}C_{il\beta}+M_{ik}C_{jl\beta})+e^{2\sigma}s_0p^2
   (C_{ij\beta}C_{lk\beta}+C_{lj\beta}C_{ik\beta}+C_{il\beta}C_{jk\beta}),
\end{eqnarray}
where $M_{ij}:=K_2h_{ij}+K_4m_im_j$.

 The result follows from (\ref{ttensor2}), Proposition \ref{c1} and the
  definition of the transformed T-tensor

 \hspace{2cm} $\overline{T}_{hijk}=\overline{L} \,\overline{C}_{hij}\overline{\mid}_k
+\overline{C}_{hij}\overline{l}_k+\overline{C}_{hik}\overline{l}_j
+\overline{C}_{hjk}\overline{l}_i+\overline{C}_{ijk}\overline{l}_h.$
\end{proof}

  The transformed  T-tensor for some
important special Finsler spaces  can be deduced from the above
result.

\begin{cor}Under a  Kropina  change,
$\overline{L}= L^2/\beta$;  L being  Reimannian, the transformed
T-tensor
 takes the form:
\begin{eqnarray}\label{kropina2}
 \nonumber \overline{T}_{lijk} &=&
  \frac{2\overline{L} }{L^2b^2}(h_{li}h_{jk}+h_{lj}h_{ik}+h_{lk}h_{ij})+
  \frac{2\overline{L}^2 }{\beta
  L^2b^2}(h_{li}m_{j}m_k+h_{lj}m_{i}m_k+h_{ij}m_{l}m_k\\
   &&+h_{jk}m_{i}m_l+h_{lk}m_{i}m_j+h_{ik}m_{j}m_l)+\frac{6\overline{L}^3 }{\beta^2L^2b^2}m_lm_im_jm_k.
\end{eqnarray}
\end{cor}
It is to be noted  that the above result is  also obtained by
Shibata \cite{r2.12}.

\begin{cor}Under a conformal change $\overline{L}=e^\sigma L$, the transformed  T-tensor
 takes the form
 \begin{eqnarray*}
  \overline{T}_{lijk} &=& e^{3\sigma}T_{lijk}
\end{eqnarray*}
\end{cor}

\begin{cor}Under a Randers change
$\overline{L}= L+\beta$, $L$ being Riemannian, the T-tensor
 takes the form:
\begin{eqnarray*}
  \overline{T}_{lijk} &=& -
  \frac{
  \Theta_1}{4L^3}(h_{li}h_{jk}+h_{lj}h_{ik}+h_{lk}h_{ij}),
\end{eqnarray*}
where $\Theta_1:=L^2b^2+\beta^2+2L\beta$.
\end{cor}
The above case has been studied by Matsumoto \cite{r2.8}.

\begin{cor}Under a $\beta$-conformal change $\overline{L}=e^\sigma L+\beta$; $L$ being Finslerian,
the transformed T-tensor
 takes the form:
 \begin{eqnarray*}
  \overline{T}_{lijk} &=&\frac{e^\sigma \overline{L}^2}{L^2}T_{lijk} -\frac{
  e^\sigma \Theta}{4L^3}(h_{li}h_{jk}+h_{lj}h_{ik}+h_{lk}h_{ij})\\
  &&+\frac{e^\sigma\overline{L}}{2L}(C_{lij}m_k+C_{ijk}m_l+C_{jlk}m_i+C_{lik}m_j) \\
   &&-\frac{e^\sigma \overline{L}}{2L}(h_{ij}C_{lk\beta}+h_{jl}C_{ik\beta}+h_{il}C_{jk\beta}
   +h_{lk}C_{ij\beta}+h_{jk}C_{il\beta}+h_{ik}C_{jl\beta}),
\end{eqnarray*}
where $\Theta:=L^2b^2+\beta^2+2e^\sigma L\beta$.
\end{cor}

\begin{cor}\label{ttensor} Under a $\beta$-conformal change, a necessary condition for the  vanishing
of  the  transformed T-tensor  is that
$$T=\frac{(n^2-1)\Theta}{4L\overline{L}^2}+\frac{(n-1)L}{\overline{L}}C_\beta,$$
where $T=g^{lj}g^{ik}T_{lijk}$.
\end{cor}

%%%%%%%%%%%%%%%%%%%%%%%%%%%%%%%%%%%%%%%%%%%%% 4. The b-condition %%%%%%%%%%%%%%%%%%%%%%%%%%%%555%%%%%%%%%%%%%%%%%%%%%%%%%%%%%

\Section{\textbf{The b-condition}} In this section we introduce
and investigate what we call the b-condition. We study the effect
of subjecting some special Finsler spaces to this condition. In
the following we assume that we are given  a generalized
$\beta$-conformal change
 $\overline{L}=f(e^\sigma L,\beta)$ with $\beta=b_iy^i=b^iy_i$.

A Finsler manifold $(M,L)$ is said to satisfy the b-condition if
$$b^iC_{ijk}=0. $$

\begin{thm}\label{b-condition1}
For $n>2$, the following two assertions are
equivalent:\begin{description}
    \item[(a)] The b-condition is invariant under a generalized
    $\beta$-conformal change.

    \item[(b)]The generalized
$\beta$-conformal change is an  energy  $\beta$-change.
 \end{description}
\end{thm}

\begin{proof}~\par
\noindent\textbf{(a)} $\Rightarrow$\textbf{(b):} Let $b^iC_{ijk}=0$.
Then, $b^i\overline{C}_{ijk}=0$ and we have, by Proposition
\ref{c1},
$$e^{\sigma}p_{-1}(m^2h_{jk}+2m_jm_k)
    +p_{02}m^2m_jm_k=0.$$
 Contracting  by $b^j$, we get $3e^{\sigma}p_{-1}=-m^2p_{02}$.
    Hence,
    $$e^{\sigma}p_{-1}(m^2h_{jk}-m_jm_k)=0,$$
    contracting again by $g^{jk}$, we get $$(n-2)p_{-1}=0.$$
    Since $n>2$, then $p_{-1}=0$ and hence the result follows   from Theorem
    \ref{energy}.

\noindent\textbf{(b)} $\Rightarrow$\textbf{(a):} Let the generalized
$\beta$-conformal change be an  energy $\beta$-change. Then, by
Theorem \ref{energy}, we obtain $\overline{C}_{ijk}=e^\sigma p\,
C_{ijk}$. Hence the result.
\end{proof}

\begin{thm}\label{b-condition2}
 Under a generalized Randers change, if $(M,L)$ satisfies the b-condition,  the generalized Randers
 space $(M,\overline{L})$ can not satisfy the b-condition.
\end{thm}

\begin{proof} Let $(M,L)$ satisfy the b-condition
$b^iC_{ijk}=0$. If $(M,\overline{L})$ satisfies the b-condition,
then $b^i\overline{C}_{ijk}=0$, and consequently,
$$\frac{1}{2L}b^i(2 \overline{L}C_{ijk}+h_{ij}m_k+h_{jk}m_i+h_{ik}m_j)=0,$$
or
$$\frac{ L }{2
\overline{L}^2}(m^2h_{jk}+2m_jm_k)=0,$$ which, by contraction  by
$g^{jk}$, yields  a contradiction: $n=-1$.
\end{proof}

\begin{thm}\label{matsumoto}Consider the generalized $\beta$-change (\ref{change}).
 In each of the following cases
\begin{description}
    \item[(a)]two-dimensional   Finsler space,

    \item[(b)]three-dimensional   Finsler space  satisfying  the
    condition $L(x,-y)=L(x,y)$,

    \item[(c)]quasi-C-reducible  space with  $b^ib^jQ_{ij}\neq 0$,

    \item[(d)] C-reducible space,

    \item[(e)] The  transformed  space  $(M,\overline{L})$ with $L$ Riemannian,
\end{description}
if the given Finsler space $(M,L)$ satisfies the b-condition, then
it is Riemannian.
\end{thm}

\begin{proof}~\par
\noindent The proof of \textbf{(a)} and \textbf{(b)} runs on in a
similar manner  as  given  in \cite{meguchi} for a concurrent vector
fields.

\noindent\textbf{(c)} Contracting  (\ref{paper2.8})  by $b^ib^j$, we
get
$$b^ib^jQ_{ij}C_k=0.$$ Hence,  $C_k=0$ for $b^ib^jQ_{ij}\neq 0$.

\noindent\textbf{(d)} Contracting  (\ref{c-reducible*})  by
$b^ib^j$,  we get
$$m^2C_k=0.$$
Consequently, $C_k=0.$

\noindent\textbf{(e)} Let $(M,\overline{L})$  be a Finsler space
with $(\alpha,\beta)$-metric, then
$$ \overline{C}_{ijk}=\frac{e^{\sigma}p_{-1}}{2}(h_{ij}m_k+h_{jk}m_i+h_{ki}m_j)
    +\frac{p_{02}}{2}m_im_jm_k.$$
  The condition  that  $b^i\overline{C}_{ijk}=0$ leads to
    $$e^{\sigma}p_{-1}(m^2h_{jk}+2m_jm_k)
    +p_{02}m^2m_jm_k=0.$$
   Contracting  by $b^j$, we get $3e^{\sigma}p_{-1}=-m^2p_{02}$.
    Hence,
    $$e^{\sigma}p_{-1}(m^2h_{jk}-m_jm_k)=0,$$
   which, by contracting by $g^{jk}$, yields  $$(n-2)p_{-1}=0.$$
    Thus, if $n=2$, the result follows by \textbf{(a)} and if
    $p_{-1}=0$, then $p_{02}=0$ and hence $\overline{C}_{ijk}=0$.
\end{proof}

\begin{thm}\label{C-change7}
 A  semi-C-reducible Finsler space  satisfying the b-condition
  is either Riemannian or C2-like.
\end{thm}

\begin{proof}
Contracting  (\ref{semi-c-red.}) by $b^jb^k$, we have $rm^2C_i=0$.
Since $m^2\neq 0$, then either  $r=0$, which implies that  the space
is $C_2$-like,
 or $C_i=0$, which implies that the  space is Riemannian.
\end{proof}

\begin{thm}\label{C-change11}
 If an  $S_3$-like  Finsler space $(M,L)$ satisfies the  b-condition, then
 its vertical curvature tensor $S_{hijk}$ vanishes.
\end{thm}
\begin{proof}
Contracting  (\ref{s3-likedef.}) by $b^l$, we get
\begin{equation}\label{new}
S(h_{ik}m_{j}-h_{ij}m_{k})=0.
\end{equation}
Again,  contracting (\ref{new}) by  $g^{ij} $, we have
$(n-2)Sm_k=0$. As $n>4$ and $m_k\neq 0$, it follows that  $S=0$ and
consequently, $S_{lijk}=0$.
\end{proof}

\begin{lem}\label{C-change1}
If a Finsler space satisfies the b-condition, then we have
$b^i|_h=0$ and, consequently,
$$C_{ijh}|_kb^h=C_{ijk}|_hb^h=0.$$
\end{lem}
\begin{proof} From the definition of vertical covariant derivative
of Cartan connection, we have
\begin{eqnarray*}
  b^i|_h&=& \dot{\partial}_h b^i+b^mC^i_{mh}
   = \dot{\partial}_h (b_jg^{ij})
   =(\dot{\partial}_h b_j)g^{ij}+b_j\dot{\partial}_h  g^{ij}
   =0
\end{eqnarray*}
and hence $C_{ijk}|_hb^h=C_{ijh}|_kb^h=(C_{ijh}b^h)|_k=0$.
\end{proof}
It is well-known  that if $(M,L)$ is Riemannian, then the T-tensor
vanishes. But the converse is not true in general. The next result
shows that the converse is true in the case where  $(M,L)$
satisfies the b-condition.
\begin{thm}\label{C-change13}
 A   Finsler space   satisfying the  b-condition is  Riemannian
  if and only if the T-tensor  $T_{hijk}$ vanishes.
\end{thm}

\begin{proof}
It is clear that if the space is Riemannian then the T-tensor
vanishes. On the other hand, if the T-tensor vanishes, then
$$LC_{hij}|_k
+C_{hij}l_k+C_{hik}l_j +C_{hjk}l_i+C_{ijk}l_h=0.$$ Contracting by
$b^i$,  using Lemma \ref{C-change1}, we have
$\frac{\beta}{L}C_{hjk}=0$. Hence $C_{hjk}=0.$
\end{proof}

\vspace{0.1cm}

Let us write
\begin{equation}\label{C-change14}
T_{ij}:=T_{ijhk}g^{hk}=LC_i|_j+l_iC_j+l_jC_i.
\end{equation}
By contracting  (\ref{C-change14}) by $b^i$, making use of Lemma
\ref{C-change13}, we have
$$T_{ij}b^i=(\beta/L)C_j.$$
Hence, we have
\begin{cor}\label{C-change13}
 A   Finsler space   satisfying the  b-condition is  Riemannian
   if and only if the  tensor $T_{ij}$  vanishes.
\end{cor}
\vspace{-.5cm}
%%%%%%%%%%%%%%%%%%%%%%%%%%%%%%%%%%%%%%%%%%%%%%%% 5. Projective Change and ... %%%%%%%%%%%%%%%%%%%%%%%%%%%%%%%%%%%%%%%%%%%%%%%

\Section{\textbf{Projective change and generalized $\beta$-conformal
change}}
 In this section we will be guided by Matsomoto \cite{proj.} and
 Shibata \cite{r2.12}.
 For two Finsler spaces $(M,L)$ and
 $(M,\overline{L})$ with the same underlying manifold
 $M$, if every geodesic on of $(M,L)$ is also a geodesic of
  $(M,\overline{L})$ and vice versa, the change $L\longrightarrow
 \overline{L}$ of  Finsler metrics is said to be projective. A
 geodesic on $(M,L)$ is characterized  by
$$\frac{dy^i}{dt}+2G^i=wy^i, \,\, \frac{dx^i}{dt}=y^i,$$
where $w=(d^2s/dt^2)/(ds/dt)^2$  and
$G^i(x,y)=\frac{1}{2}\gamma^i_{jk}y^jy^k$ is the canonical spray of
$(M,L)$. We are going to find out a
condition for a generalized $\beta$-conformal change to be
projective.

 Consider the left hand side of   Euler-Lagrange  equations
\begin{equation}\label{euler1}
\mathcal{E}_i:=\partial_i  L- \frac{d}{dt}(\dot{\partial}_iL)
\end{equation}
\begin{prop}\label{euler2}
 Under a generalized     $\beta$-conformal change $\overline{L}=f(e^\sigma
 L,\beta)$, the functions (\ref{euler1}) are transformed according  to
 \begin{equation}\label{phi}
 f\overline{\mathcal{E}}_i=Le^\sigma p\mathcal{E}_i+Lq_0 m_i m^r \mathcal{E}_r +\varphi_i,
 \end{equation}
\noindent where
 \begin{equation}\label{paper2.2}
 \varphi_i:=L^2e^\sigma p \,\sigma_i -(pLe^\sigma l_i-q_0\beta
 \,m_i)\sigma_0+2qF_{0i}-q_0E_{00}m_i,
 \end{equation}
$$\sigma_i:=\partial_i\sigma,  \quad
F_{ij}:=\frac{1}{2}(b_{i|j}-b_{j|i}), \quad
 E_{ij}:=\frac{1}{2}(b_{i|j}+b_{j|i}),$$
 $$\sigma_0=\sigma_iy^i,\quad F_{0i}=F_{li}y^l,\quad E_{00}=E_{ij}y^iy^j.$$
\end{prop}
\begin{proof} Making use  of the homogeneity of $f$,  $\overline{E_i}$ can be computed    as
follows.
\begin{eqnarray}\label{russion}
 \nonumber     \overline{\mathcal{E}_i}&=&\partial_i  f- \frac{d}{dt}(\dot{\partial}_if)\\
  \nonumber    &=&f_1(\sigma_i e^\sigma L+e^\sigma
      \partial_iL)+f_2(N^r_ib_r+b_{j\mid i}\,y^j)-\frac{d}{dt}(e^\sigma f_1
      l_i+f_2b_i)\\
   \nonumber   &=&f_1\sigma_i e^\sigma L+f_1 e^\sigma\partial_iL+f_2
      N^r_ib_r+f_2b_{j\mid i}\,y^j-f_1l_i e^\sigma \sigma_r y^r-f_1e^\sigma\frac{d\,l_i }{dt}
      -e^\sigma l_i\frac{d \,f_1}{dt}\\
  \nonumber    &&-b_i\frac{d \,f_2}{dt}-f_2(b_{i\mid j}\,y^j+N^r_ib_r)\\
      &=&f_1 e^\sigma \mathcal{E}_i+f_1L\sigma_i e^\sigma -f_1\sigma_0 e^\sigma
      l_i+2f_2F_{0i}-\frac{d\,f_2}{dt}m_i.
      \end{eqnarray}
Using the relation
$\displaystyle{\frac{d\,y_r}{dt}=y^s\partial_sy_r+g_{rs}
\frac{d\,y^s}{dt}}$, the last  term
$\displaystyle{\frac{d\,f_2}{dt}}$ of (\ref{russion}) is given by
      \begin{eqnarray}\label{russion2}
   \nonumber   \frac{d\,f_2}{dt}&=&f_{21}\frac{d\,\tilde{L}}{dt}+f_{22}\frac{d\,\beta}{dt}\\
      \nonumber  &=&-\frac{\beta f_{22}}{e^\sigma L}(e^\sigma y^r \partial_r L+e^\sigma  l_r\frac{d\,y^r}{dt}+L\frac{d\,e^\sigma}{dt})
      +f_{22}(E_{00}+N^r_sb_ry^s+b_r\frac{d\,y^r}{dt})\\
      &=&f_{22}E_{00}-Lf_{22}m_i m^r\mathcal{E}_r -\beta f_{22}\sigma_0.
\end{eqnarray}
      Now, substituting  (\ref{russion2}) into (\ref{russion}), we
get\\
     \emph{$ f\overline{\mathcal{E}_i}=Le^\sigma p\mathcal{E}_i+q_0 L m_i m^r\mathcal{E}_r+L^2e^\sigma p \,\sigma_i -(pLe^\sigma l_i-q_0\beta
      \,m_i)\sigma_0+2qF_{0i}-q_0E_{00}m_i$}\\
    \emph{${ \hspace{.6cm}} =L p \,e^\sigma \mathcal{E}_i+q_0Lm_i m^r\mathcal{E}_r+\varphi_i.$}
\end{proof}

\begin{thm}\label{projective} A
generalized     $\beta$-conformal change is projective if and only
if the vector $\varphi_i$  vanishes.
\end{thm}
\begin{proof} Let the generalized $\beta$-conformal change be
projective.  Then,   $\mathcal{E}_i=0$ is equivalent to
$\overline{\mathcal{E}}_i=0$ and consequently, $\varphi_i=0$ by
(\ref{phi}).

Conversely, if $\varphi_i=0$, then (\ref{phi}) shows that
$\mathcal{E}_i=0$ implies  $\overline{\mathcal{E}}_i=0$. On the
other hand, if $\overline{\mathcal{E}}_i=0$ and $\varphi_i=0$, then
$e^\sigma p   \mathcal{E}_i+q_0 m_i m^r\mathcal{E}_r=0$. Contracting
the last equation by $m^i$, taking into account that  $e^\sigma
p+m^2 q_0\neq0$, we get $\mathcal{E}_rm^r =0.$ Consequently,
$\mathcal{E}_i =0$.
\end{proof}

From the above theorem,  we retrieve   the    following two
results due to Shibata \cite{r2.13} and Hashiguchi and Ichijo
\cite{hashiguchi} respectively.

\begin{cor}
A $\beta$-change is projective if and only if \quad
$2qF_{0i}=q_0E_{00}m_i$.
\end{cor}
\begin{cor}\label{hashiguchi}
A Randers change is projective if and only if $F_{0i}=0$, that is,
$b_i$ is gradient.
\end{cor}

The following two results are a generalized version of Shibata's
result \cite{r2.13} and Matsumoto's result  \cite{proj.}.

\begin{thm}\label{paper2.5}
Assume that the generalized $\beta$-conformal change (\ref{change})
is projective and $L$ is  Minkowskian, then the Weyl  torsion
$\overline{W}^h_{ij}$ and the  Douglas  tensor
$\overline{D}^h_{ijk}$ of $(M,\overline{L})$ vanish. Consequently,
$(M,L)$ with dim  $M>2$ is projectively flat.
\end{thm}

\begin{proof} The  Weyl torsion tensor  is given by \cite{proj.}:
$$W^h_{ij}=\overcirc{R}^h_{ij}+\frac{1}{n+1}\mathfrak{A}_{(i,j)}\{y^h\,\overcirc{R}_{ij}
+\delta^h_i\,\overcirc{R}_{j} \},$$ where
  \, $\overcirc{R}_{ij}=\overcirc{R}^h_{ijh}$,\,
$\overcirc{R}_{j}=\frac{1}{n+1}(n\,\,\overcirc{R}_{0j}+\overcirc{R}_{j0})$
 and \, $\overcirc{R}^h_{ijk}$ is the h-curvature of the Berwald
connection. Since $(M,L)$ is Minkowskian, then
$\overcirc{R}^h_{ijk}=0$, and so
 \, $\overcirc{R}_{ij}=\overcirc{R}_{i}=0$. Consequently,
$W^h_{ij}=0$. By the invariance of $W^h_{ij}$ under a projective
change, we have $\overline{W}^h_{ij}=0$.

The   Douglas tensor is given by \cite{proj.}:
$$D^h_{ijk}=\overcirc{P}^h_{ijk}+\frac{1}{n+1}(y^h\,\overcirc{P}_{ij}{\stackrel{\circ}|_k}
+\mathfrak{S}_{(i,j,k)}\{\delta^h_i\,\overcirc{P}_{jk} \}),$$
where\, $\overcirc{P}_{ij}=\overcirc{P}^h_{ijh}$,\,\,
$\overcirc{P}^h_{ijk}$ is the hv-curvature of the Berwald connection
and $\stackrel{\circ}|$ denotes the vertical covariant derivative
with respect to the Berwald connection $G^h_{ij}$.
 Since $(M,L)$ is
Minkowskian, then \, $\overcirc{P}^h_{ijk}=0$, and so\,
 $\overcirc{P}_{ij}=0$. Consequently, $D^h_{ijk}=0$. By the
invariance of $D^h_{ijk}$ under a projective change, we have
$\overline{D}^h_{ijk}=0$.

Finally, as
  $W^h_{ij}=0$, $D^h_{ijk}=0$ and dim $M>2$, $(M,L)$
  is thus projectively flat \cite{proj.}.
\end{proof}

\begin{thm}\label{paper2.6}
Assume that the generalized $\beta$-conformal change is projective
and $L$ is  Riemannian,  then the projective
 hv-curvature tensor $\overline{D}^h_{ijk}$ of
 $(M,\overline{L})$ vanishes.
\end{thm}

\begin{proof}
Since $(M,L)$ is Riemannian, then \, $\overcirc{P}^h_{ijk}=0$, and
$\overcirc{P}_{ij}=0$. Consequently, $D^h_{ijk}=0$. By the
invariance of $D^h_{ijk}$ under a projective change, we have
$\overline{D}^h_{ijk}=0$.
\end{proof}

\begin{thm}\label{proj.2}
If $\varphi_i=0$, then $(M,\overline{L})$ is of scalar curvature
 if and only if $(M,L)$  is
  of scalar curvature.
\end{thm}

\begin{proof}
According  to Szab\'{o} \cite{szabo}, a Finsler space is of scalar
curvature if and only if $W^h_{ij}=0$ vanishes identically. Let
$\varphi_i=0$, then by Theorem \ref{projective} the  generalized
$\beta$-conformal change is projective. Now, let $(M,L)$  be of
scalar curvature, then  $W^h_{ij}=~0$. But
$\overline{W}^h_{ij}=W^h_{ij}$, hence, $\overline{W}^h_{ij}=0$.
Consequently, $(M,\overline{L})$ is of scalar curvature. Conversely,
let $(M,\overline{L})$ be  of scalar curvature, then
$\overline{W}^h_{ij}=0$ which leads to $W^h_{ij}=0$, hence,
$({M},L)$ is of scalar curvature.
\end{proof}
In the Riemannian case the term \lq\lq of scalar curvature\rq\rq
  \, reduces to the term \lq\lq of constant curvature\rq\rq. Thus , we
generalize Yasuda and Shimada's
 result \cite{r2.17}.
\begin{cor}\label{paper2.4}
Under a generalized $\beta$-conformal change, if  $\varphi_i=0$ and
$(M,L)$  is  Riemannian, then the Finsler space $(M,\overline{L})$
is  of scalar curvature if and only if $(M,L)$ is of constant
curvature.
\end{cor}

\vspace{-2cm} \providecommand{\bysame}{\leavevmode\hbox
to3em{\hrulefill}\thinspace}
\providecommand{\MR}{\relax\ifhmode\unskip\space\fi MR }
\providecommand{\MRhref}[2]{%
  \href{http://www.ams.org/mathscinet-getitem?mr=#1}{#2}
} \providecommand{\href}[2]{#2}

\end{document}